\documentclass[twoside]{article}
\usepackage{amsmath,amsthm,amsfonts,amssymb}
\def\eqnreset{\setcounter{equation}{0}}

\newtheorem{lemma}{Lemma}[section]
\newtheorem{thm}[lemma]{Theorem}
\newtheorem{defi}[lemma]{Definition}
\newtheorem{cor}[lemma]{Corollary}
\newtheorem{prop}[lemma]{Proposition}

\def\medno{\medbreak\noindent}

\def\qed{~\hfill$\square$\medbreak}

\def\naam#1{\label{#1}}
\def\refer#1{\ref{#1}}

\def\bib#1{\cite{#1}}

\def\text#1{\;\;\;\;{\rm \hbox{#1}}\;\;\;\;}
\def\qquad{\quad\quad}

\def\minspace{\vspace{-2mm}}

\def\itema{\vspace{-1mm}\item[{\rm (a)}]}
\def\itemb{\minspace\item[{\rm (b)}]}

\def\msy#1{{\mathbb #1}}
\def\C{{\msy C}}
\def\N{{\msy N}}
\def\Z{{\msy Z}}
\def\R{{\msy R}}
\def\D{{\msy D}}

\def\ga{\alpha}

\def\ge{\varepsilon}
\def\gf{\varphi}

\def\gl{\lambda}

\def\gs{\sigma}

\def\gD{\Delta}

\def\gS{\Sigma}

\def\got#1{\mathfrak #1}
\def\fa{{\got a}}

\def\fg{{\got g}}
\def\fh{{\got h}}

\def\fk{{\got k}}

\def\fn{{\got n}}

\def\fp{{\got p}}
\def\fq{{\got q}}

\def\to{\rightarrow}
\def\Re{{\rm Re}\,}

\def\inp#1#2{\langle#1\,,\,#2\rangle}

\def\Ad{{\rm Ad}}
\def\End{{\rm End}}
\def\Hom{{\rm Hom}}

\def\ad{{\rm ad}}
\def\after{\,{\scriptstyle\circ}\,}

\def\iq{{\rm q}}

\def\iC{{\scriptscriptstyle \C}}
\def\iR{{\scriptscriptstyle \R}}

\def\cA{{\cal A}}

\def\cC{{\cal C}}

\def\cF{{\cal F}}
\def\cH{{\cal H}}

\def\cM{{\cal M}}

\def\cO{{\cal O}}
\def\cP{{\cal P}}
\def\cR{{\cal R}}
\def\cS{{\cal S}}
\def\cT{{\cal T}}

\def\cW{{\cal W}}

\mathcode`:="603A
\def\col{\,:\,}

\def\spX{{\rm X}}
\def\Ci{C^\infty}
\def\Vtau{V_\tau}
\def\Cartan{\theta}
\def\Cci{C^\infty_c}
\def\faq{\fa_\iq}
\def\faqdc{\fa_{\iq\iC}^*}
\def\NKaq{N_K(\faq)}
\def\ZKaq{Z_K(\faq)}
\def\faFq{\fa_{F\iq}}
\def\Ind{{\rm Ind}}
\def\Pmin{P_\emptyset}
\def\oC{{}^\circ \cC}
\def\WKH{W_{K\cap H}}
\def\Aq{A_\iq}
\def\Vtauempty{V_\tau^\emptyset}
\def\nC{C^\circ}
\def\dotvar{\,\cdot\,}
\def\rmI{{\rm I}}
\def\dE{E^*}
\def\Mer{{\mathcal M}}
\def\Fou{{\mathcal F}}
\def\Foumc{\cF_\emptyset}
\def\faqd{\fa_\iq^*}
\def\Wave{{\mathcal J}}
\def\Lau{{\mathcal L}}
\def\ev{{\rm ev}}
\def\laur{{\rm laur}}
\def\supp{{\rm supp}\,}
\def\Wavemc{\Wave_\emptyset}
\def\spX{{\rm X}}
\def\spXp{{\spX_+}}
\def\Eps{{E_{+,s}}}
\def\Epone{{E_+}}
\def\faFqdcperp{{}^*\fa_{F\iq\iC}^*}
\def\staFqd{{}^*\fa_{F\iq}^*}
\def\faFqdp{\fa_{F \iq}^{*+}}
\def\faFqdc{\fa_{F\iq\iC}^*}
\def\staFq{{}^*\fa_{F\iq}}
\def\staFqdc{{}^*\fa_{F\iq\iC}^*}
\def\faFqd{\fa_{F\iq}^*}
\def\DX{{\mathbb D}(\spX)}
\def\AFqp{A_{F\iq}^+}
\def\spXF{\spX_F}
\def\Aqp{{A_\iq^+}}
\def\nE{E^\circ}
\def\Vtau{V_\tau}
\def\oCtau{{}^\circ\cC(\tau)}
\def\AC{{\rm AC}\,}
\def\Hyp{{\mathcal H}}
\def\PW{{\rm PW}\,}
\def\PWM{{\rm PW}_{\! M}}
\def\CiM{C^\infty_M}
\def\start{{{}^*t}}
\def\spXFp{\spX_{F+}}
\def\bp{{}^\backprime}
\def\diag{{\rm diag}\,}
\def\faqd{\fa_\iq^*}
\def\spXmin{\spX_\emptyset}
\def\taumin{\tau_\emptyset}
\def\faqp{\fa_\iq^+}
\def\Aqp{A_\iq^+}
\def\parabs{\cP_\gs}
\def\dEps{E^*_{+,s}}
\def\faqdp{\fa_\iq^{*+}}
\begin{document}
\title{Eisenstein integrals and induction of relations}
\author{E.\ P. van den Ban}\date{}
\maketitle
\footnote{{\bf 2000 Mathematics Subject Classification:} 
Primary 22E30; Secondary 43A85, 22E45.}
\footnote{Keywords: Symmetric space, Eisenstein integral, induction,
Paley-Wiener theorem, Arthur-Campoli relations, Plancherel theorem.}
\pagestyle{myheadings}
\markboth{E.\ P.\ van den Ban}{Eisenstein integrals and induction of relations}
\vspace{-5mm}
\mbox{}\newline\centerline{In Honor of Jacques Carmona}
\vspace{5mm}
\begin{abstract}
In this article I will give a survey of joint work with Henrik
Schlichtkrull on the induction of certain relations among
(partial) Eisenstein integrals for the minimal principal series of
a reductive symmetric space. I will discuss the application of
this principle of induction to the proof of the Fourier inversion
formula in \bib{BSfi} and to the proof of the Paley-Wiener theorem
in  \bib{BSpw}. Finally, the relation with
the Plancherel decomposition will be discussed.
\end{abstract}
\section{Introduction}
\eqnreset
Let $\spX = G/H$ be a reductive symmetric space, with $G$ a real
reductive group of Harish-Chandra's class and $H$ an open subgroup of the
group $G^\gs$ of fixed points for an involution $\gs$ of $G.$ Thus,
$G^\gs_e \subset H \subset G^\gs,$ with
$G^\ga_e$ the identity component of $G^\gs.$

There exists a Cartan involution $\Cartan$ of $G$ that
commutes with $\gs.$ The associated
maximal compact subgroup $K:= G^\Cartan$ is invariant under $\gs.$

There are two important classes of examples of reductive symmetric spaces.
The first class, with $H$ compact, consists of the Riemannian symmetric spaces.
Here we take $\theta = \gs$ and $K = H.$ The second consists of the real reductive
groups of Harish-Chandra's class. Given such a group $\bp G,$ let $G = \bp G \times \bp G,$
let $\gs: G \to G,$ $(x,y)\mapsto (y,x),$ and let $H = G^\gs = \diag(\bp G).$ Then $\spX$
equals $\bp G,$ equipped with the left times right action of $G.$
We may take $\Cartan = \bp \Cartan \times \bp \Cartan,$ with $\bp \Cartan$ a Cartan involution of
$\bp G.$ Accordingly, $K = \bp K \times \bp K,$ with $\bp K$ maximal compact in $\bp G.$

We are interested
in the analysis of $K$-finite functions on $\spX.$ For this it is convenient
to fix a finite dimensional unitary representation $(\tau, \Vtau)$ of $K$ and to
consider the space
\begin{equation}
\naam{e: smooth spherical functions}
\Ci(\spX\col \tau) : = [\Ci(\spX) \otimes \Vtau]^K
\end{equation}
of smooth $\tau$-spherical functions on $\spX.$
Alternatively, we view $\Ci(\spX\col \tau)$
as the space of smooth functions $f: \spX \to \Vtau$ transforming according to the rule
$f(kx) = \tau(k) f(x),$ for $x \in \spX$ and $k \in K.$ The subspace of compactly supported
functions in (\refer{e: smooth spherical functions}) is denoted by $\Cci(\spX\col \tau).$

As usual, we denote Lie groups by Roman capitals,
and the associated
Lie algebras by the corresponding German lowercase letters. The involutions
$\gs$ and $\Cartan$ of $G$ give rise to involutions of the Lie algebra $\fg,$ which
are denoted by the same symbols.
Accordingly, we write
$$
\fg = \fk \oplus \fp = \fh \oplus \fq
$$
for the decompositions of $\fg$ into the $+1$ and $-1$ eigenspaces for $\Cartan$ and
$\gs,$ respectively. Let $\faq$ be a maximal abelian subspace of $\fp \cap \fq.$
The set $\gS = \gS(\fg, \faq)$ of restricted roots of $\faq$ in $\fg$ is a (possibly non-reduced)
root system. Let $\gS^+$ be a positive system, $\faqp$ the associated positive
chamber, $\Aqp: = \exp \faqp,$ and $\gD$ the associated collection of
simple roots. The Weyl group $W$ of $\gS$ is canonically isomorphic with $\NKaq/\ZKaq.$
Each subset $F \subset \gD$ determines a standard parabolic subgroup $P_F$ of $G$
as follows.
Let $\faFq$ be the intersection of the root hyperplanes $\ker \ga,$ for $\ga \in F,$
and let $M_{1F}$ be the centralizer of $\faFq$ in $G.$ Moreover, let
$$
\fn_F = \bigoplus_{\ga \in \gS^+\setminus \Z F} \fg_\ga,\text{and} N_F: = \exp \fn_F.
$$
Then $P_F = M_{1F} N_F.$ Let $M_{1F} = M_F A_F$ according to the Langlands decomposition of $P_F,$
then
$\faFq$ is the intersection of $\fa_F,$ the Lie algebra of $A_F,$ with $\fq.$ In particular,
$\faq$ is the intersection of $\fa := \fa_\emptyset$ with $\fq.$ The group $M_F$
is a reductive group of Harish-Chandra's class; accordingly, the homogeneous
space $\spX_F: = M_F / M_F \cap H$ belongs to the class of symmetric spaces
considered.

Since $\gs$ and $\Cartan$ commute, the composition $\gs \Cartan$ is an involution of $G.$
Its derivative
restricts to the identity
on $\faq;$ therefore, the involution $\gs\Cartan$ leaves each of the standard parabolic
subgroups $P_F$ invariant. Let $\parabs$ denote the collection of $\gs\Cartan$-stable
parabolic subgroups of $G$ containing $\Aq.$ Then $W$ acts on $\parabs$
in a natural way. Each element of $\parabs$ is $W$-conjugate to a unique standard
parabolic subgroup $P_F.$ Finally, each minimal element of $\parabs$
is $W$-conjugate to  $P_\emptyset.$

In this article we will discuss relations of a certain type
between (partial) normalized Eisenstein integrals for $\Pmin.$ These Eisenstein integrals,
denoted $\nE(\gl \col \dotvar),$ are essentially
sums of matrix coefficients of induced representations of the form
$\Ind_{\Pmin}^G(\xi \otimes \gl \otimes 1),$ with $\xi$ an irreducible finite dimensional unitary
representation of $M_\emptyset$ and with $\gl \in \faqdc.$ These
induced representations form
the minimal principal series of $\spX.$  Induction
of relations describes how relations of a certain type between the Eisenstein integrals
$\nE(\gl)$ (or more generally between partial Eisenstein integrals)
are induced by similar relations between the similar integrals for $\spXF.$

In terms of the mentioned Eisenstein
integrals we define a Fourier transform $\Foumc.$ 
Applied to a function $f \in \Cci(\spX\col \tau)$
the Fourier transform gives an element of $\cM(\faqdc)\otimes \oC,$
where $\cM(\faqdc)$ denotes the space of meromorphic functions on $\faqdc$ and where
$\oC = \oCtau$  is a certain finite dimensional Hilbert space. It is a main
result of \bib{BSmc} that the Fourier transform $\Foumc$ is injective on $\Cci(\spX\col \tau).$
Accordingly, two natural problems arise.
\begin{enumerate}
\itema
To retrieve
$f$ from its Fourier transform $\Foumc f$; this is the problem of {\bf Fourier inversion}.
\itemb
To characterize the image of
$\Foumc(\Cci(\spX\col \tau))$ in a way that generalizes
the {\bf Paley-Wiener theorem} of J.\ Arthur, \bib{Arthur}.
\end{enumerate}
In the answers to these related questions,
given in \bib{BSfi} and \bib{BSpw}, respectively,
the principle of induction of relations plays a fundamental role.

\section{Eisenstein integrals}
\naam{s: Eis}
\eqnreset
As said, Eisenstein integrals for $\Pmin$ are essentially sums of $K$-finite matrix coefficients
of principal series representations of the form $\Ind_{\Pmin}^G(\xi\otimes \gl \otimes 1).$
We will now give their precise definition. To keep the exposition as light as possible,
we make the following
\medno
{\bf Simplifying assumption\ }{\ }\sl The manifold $G/H$ has precisely one open $\Pmin$-orbit.
\rm
\medno
This assumption
is only made for purposes of exposition, it is not
necessary for the development of the theory. In the general
situation, there are finitely many open $\Pmin$-orbits, naturally
parametrized by $W/\WKH,$ where $\WKH$ denotes the subgroup of $W$ consisting
of elements that are contained in the natural image of $\NKaq \cap H.$
The simplifying assumption is satisfied in the Riemannian case as well as in case of the group.

We define $\oC = \oC(\tau)$ by
$$
\oC := \Ci(\spXmin\col \taumin),
$$
the space of smooth $\taumin$-spherical functions $\spXmin \to \Vtau;$
here $\taumin: = \tau|_{K \cap M_\emptyset}.$ By compactness of
$\spXmin,$ it follows that $\oC$ is finite dimensional. Moreover,
$\oC = L^2(\spXmin \col \taumin).$

Given $F \subset \gD$ we define $\rho_F \in \faFqd$ by $\rho_F = 
\frac12 {\rm tr}\,(\ad(\dotvar)|_{\fn_F}).$
In particular, we put $\rho = \rho_\emptyset.$
Let $\psi \in \oC$ and $\gl \in \faqdc.$ We define the function
$\psi_\gl: G \to \Vtau$ by
\begin{eqnarray*}
\psi_{\gl}(x) &=& a^{\gl + \rho} \psi(m)\quad\text{for} \quad x \in manH,\;\;\;\;(m,a,n) \in
M_\emptyset \times A_\emptyset \times N_\emptyset,\\
&=&0\qquad\;\;\;\;\;\;\;\;\;\;\,\text{for} \;\;\;\,x \in G\setminus \Pmin H.
\end{eqnarray*}
We equip $\fg$ with a non-degenerate $\Ad(G)$-invariant bilinear form $B$
that is negative definite on $\fk$ and positive definite on $\fp$ and
for which $\fh$ and $\fq$ are orthogonal. Then $B$ induces a positive definite
inner product $\inp{\dotvar}{\dotvar}$ on $\faqd$ which is extended to a complex
bilinear form on $\faqdc.$
For $R\in \R,$ we define
\begin{equation}
\naam{e: faqd Pmin R}
\faqd(\Pmin, R):= \{ \nu \in \faqdc \mid \inp{\Re \nu}{\ga} < R, \;\forall \, \ga \in \Sigma^+ \}.
\end{equation}
For $\gl \in    - \rho +  \faqdc(\Pmin, 0) $
we define the {\bf Eisenstein integral}
$E(\psi \col \gl\col \dotvar),$ also denoted 
$E(\Pmin \col \psi \col \gl\col \dotvar),$ by
\begin{equation}
\naam{e: defi Eisenstein}
E(\psi \col \gl\col x) = \int_K \tau(k) \psi_\gl(k^{-1}x)\; dk,
\end{equation}
for $x \in \spX.$ The following result is due to \bib{Bps2}, Prop.\ 10.3.

\begin{prop}
The integral (\refer{e: defi Eisenstein})
is absolutely convergent for $\gl \in - \rho + \faqd(\Pmin,0)$
and defines a holomorphic function of $\gl$ with values in
$\Ci(\spX\col \tau).$  Moreover, it extends to a meromorphic function of
$\gl \in \faqdc$
with values in the space $\Ci(\spX\col \tau).$ The singular locus
of this meromorphic extension is a locally finite union
of hyperplanes of the form $\gl_0 + (\ga^\perp)_\iC,$
with $\gl_0 \in \faqdc$ and $\ga \in \gS.$
\end{prop}

For generic $\gl \in \faqdc$ with $\inp{\Re \gl - \rho}{\ga}>0$ for all $\ga \in \gS^+$
we have that
$$
\lim_{\buildrel {a \to \infty}\over {a\in \Aqp}}
a^{-\gl + \rho} E(\psi\col \gl \col a) = [C(1: \gl) \psi](e)
$$
with $C(1 \col \gl) = C_{\Pmin|\Pmin}(1 \col \gl)  \in \End(\oC)$ a meromorphic function of
$\gl$ that extends meromorphically to all of $\faqdc;$ see \bib{Bps2}, Sect.\ 14.
Since the function $\gl \mapsto \det C(1\col \gl)$
is not identically zero, we may define the {\bf normalized Eisenstein
integral}
$$
\nE(\psi\col \gl\col x): = E(C(1\col \gl)^{-1}\psi \col \gl \col x),
$$
see \bib{Bps2}, Sect.\ 16, and \bib{BSft}, Sect.\ 5, for details. The definition generalizes
that of Harish-Chandra, \bib{HCeis}, p.\ 135, in the case of the group.

The normalized Eisenstein integral
$\nE(\psi\col \gl)$ is a meromorphic function of  $\gl \in \faqdc$ with
values in $\Ci(\spX\col \tau).$ Its asymptotic behavior
is described by the following theorem.
Given $a \in \Aq$ we write $z(a)$ for the point in $\C^\gD$ with components
$z(a)_\ga = a^{-\ga},$ for $\ga \in \gD.$
Let $D \subset \C$ denote the complex unit disc. Then $z$ maps $\Aqp$ into $D^\gD.$
If $\Omega$ is a complex analytic manifold, then by $\cO(\Omega)$ we denote the algebra of
holomorphic functions on $\Omega.$
Let $\Vtauempty$ denote the space of $M_\emptyset \cap K \cap H$-fixed elements in $\Vtau.$

\begin{prop}
\naam{p: asymptotics Eis}
There exists a unique meromorphic function $\gl \mapsto \Phi_\gl$ with values
in $\cO(D^\gD)\otimes \End( \Vtauempty)$ and, for $s \in W,$
unique meromorphic $\End(\oC)$-valued
meromorphic functions $\faqdc \ni \gl \mapsto \nC(s \col \gl)$ such that, for
all $\psi \in \oC,$
\begin{equation}
\naam{e: expansion nE}
\nE(\psi\col\gl\col a)  =\sum_{s \in W} a^{s\gl - \rho}
 \Phi_{s\gl}(z(a)) \; [\nC(s \col \gl)\psi](e),
\end{equation}
for $a \in \Aqp,$ as a meromorphic identity in the variable $\gl \in \faqdc.$ The meromorphic
functions $\gl \mapsto \Phi_\gl$ and $\gl \mapsto \nC(s\col \gl),$ for $s \in W,$
all have a singular locus that is a locally finite union of
hyperplanes  of the form
$\gl_0 + (\ga^\perp)_\iC,$ with $\gl_0 \in \faqdc$ and $\ga \in \gS.$
\end{prop}

For a proof of this result, we refer the reader to
\bib{BSexp}, Sect.\ 11, and \bib{BSanfam}, Sect.\ 14.
From the above result it follows in particular that $\gl \mapsto \nE(\psi\col \gl)$
is a meromorphic $\Ci(\spX\col \tau)$-valued function with
singularities along hyperplanes of the form $\gl_0 + (\ga^\perp)_\iC,$ with
$\gl_0 \in \faqdc$ and $\ga \in \gS.$

We note that it follows from the definition of the normalized Eisenstein integral that
$$
\nC(1\col \gl) = \rmI_{\,\oC},
$$
as a meromorphic
identity in the variable $\gl \in \faqdc.$
Besides in the expansion (\refer{e: expansion nE})
the normalized $c$-functions $\nC(s\col\dotvar)$ also appear in
the following functional equation for the Eisenstein integral
\begin{equation}
\naam{e: functional equation}
\nE(s\gl \col x) \nC(s\col \gl) = \nE(\gl \col x),
\end{equation}
for every $x \in \spX,$ as an identity of meromorphic functions in
the variable $\gl \in \faqdc.$

The following result is crucial for the further development of the theory.
\begin{thm}
\naam{t: Maass-Selberg}
{\bf (Maass-Selberg relations)\ }{\ } For each $s \in W,$
$$
\nC(s \col - \bar \gl)^* \nC(s \col \gl) = \rmI_{\oC},
$$
as a meromorphic identity in the variable $\gl \in \faqdc.$
\end{thm}

In the case of the group, the terminology Maass-Selberg relations was introduced by Harish-Chandra,
because of striking analogies with the theory of automorphic forms.
In the mentioned setting of the group Harish-Chandra derived
the relations for the $c$-functions associated with arbitrary parabolic subgroups, see
\bib{HC3}.
In the present setting the above result is due to \bib{Bps2}, Thm.\ 16.3,
see also \bib{Bsd}.
The result has been generalized to $c$-functions associated with arbitrary
$\gs\Cartan$-stable parabolic subgroups by P.~Delorme \bib{Dtr}, see
also \bib{CDn}.
It plays a crucial role in Delorme's proof of the Plancherel formula, see \bib{Dpl},
as well in the proof of the Plancherel formula by myself and Schlichtkrull, see
\bib{BSpl1} and \bib{BSpl2}. Recently the last mentioned authors have been
able to obtain the Maass-Selberg relations for arbitrary parabolic subgroups
from those for the minimal one, see \bib{BSpl1}. The proof in the latter paper
is thus independent from the one by Delorme.

From the Maass-Selberg relations, combined with
the information that  the singular locus of the meromorphic $c$-functions
is a locally finite union of translates of root hyperplanes,
the following result is an easy consequence.

\begin{cor}
Let $s \in W.$ The normalized $c$-function $\nC(s\col \dotvar)$ is regular on $i\faqd.$ Moreover,
for $\gl \in i\faqd,$ the endomorphism $\nC(s \col \gl) \in \End(\oC)$ is unitary.
\end{cor}

From this result  and an asymptotic analysis involving induction with respect
to the split rank of $\spX,$ i.e., $\dim \faq,$
it can be shown that the normalized Eisenstein
integrals are regular for imaginary values of $\gl.$ This is the main motivation
for their definition.

\begin{thm}\naam{t: regularity}
{\bf (Regularity theorem)\ }{\ }
Let $\psi \in \oC.$ The Eisenstein integral $\nE(\psi \col \gl)$ is meromorphic in $\gl \in \faqdc$
with a singular locus  disjoint from $i\faqd.$
\end{thm}

The above result is due to \bib{BSft}, p.\ 537, Thm.\ 2.
A different proof of the regularity
theorem has been given by \bib{BCDn}. The latter approach was generalized
to arbitrary $\gs\Cartan$-stable parabolic subgroups by J.\ Carmona and P.\ Delorme,
yielding the regularity theorem for Eisenstein integrals as a consequence
of the Maass-Selberg relations in that setting;
see \bib{CDn}, Thm.\ 3 (i).

\section{Fourier inversion}
\eqnreset
The regularity theorem allows us to define a Fourier transform that is regular
for imaginary values of the spectral parameter $\gl.$ For its definition
it is convenient to define
$\nE(\gl \col x) \in \Hom(\oC, \Vtau)$ by
$$
\nE(\gl \col x) \psi := \nE(\psi\col \gl \col x).
$$
In addition, we define the dualized Eisenstein integral by conjugation,
\begin{equation}
\naam{e: defi dual Eis}
\dE(\gl \col x) := \nE(-\bar \gl \col x)^* \in \Hom(\Vtau, \oC),
\end{equation}
for $x \in \spX,$ as a meromorphic function of $\gl \in \faqdc.$ We now define
the (most-continuous) {\bf Fourier transform} $\Foumc f$ of a function
$f \in \Cci(\spX\col \tau)$ to be the meromorphic function in $\Mer(\faqdc) \otimes \oC$
given by
\begin{equation}
\naam{e: defi Foumc}
\Foumc f(\gl) := \int_\spX\; \dE(\gl \col x) f(x)\; dx,\qquad (\gl \in \faqdc).
\end{equation}
It follows from (\refer{e: functional equation})
combined with the definition of $\dE(\gl\col x)$
and the Maass-Selberg relations that, for each $s \in W,$
\begin{equation}
\naam{e: transformation property Fou}
\Foumc f(s\gl) =  \nC(s\col \gl) \Foumc f(\gl).
\end{equation}
It follows from the regularity theorem that the Fourier transform $\Foumc f$ is a regular
function on $i\faqd.$ The following theorem is one of the main results
of \bib{BSmc}, see loc.\ cit., Thm.\ 15.1.

\begin{thm}
\naam{t: injectivity Foumc}
The Fourier transform $\Foumc$ is injective on $\Cci(\spX\col \tau).$
\end{thm}

There exists a notion of {\bf Schwartz space} $\,\cC(\spX\col \tau),$
which is the proper generalization of Harish-Chandra's Schwartz space for the group,
see \bib{Bps2}, Sect.\ 17. It has the property that $\Foumc$
extends to a continuous linear map from $\cC(\spX\col \tau)$ into the Euclidean
Schwartz space $\cS(i\faqd) \otimes \oC,$
see \bib{BSft}, p.\ 573, Cor.\ 4.
We emphasize that the extended Fourier transform is in general
not injective on the Schwartz space.
More precisely,
there is a continuous linear {\bf wave packet transform} $\Wavemc:\; \cS(i\faqd) \otimes \oC
\to \cC(\spX\col \tau),$ defined by the formula
\begin{equation}
\naam{e: defi wave packet}
\Wavemc \gf(x) = \int_{i\faqd} \; \nE(\gl\col x) \gf(\gl)\; d\gl, \qquad (x \in \spX),
\end{equation}
for $\gf \in \cS(i\faqd) \otimes \oC,$ see \bib{BCDn}, Thm.\ 1.
Here $d\gl$ denotes Lebesgue measure on $i\faqd,$ suitably normalized.
Furthermore, in  \bib{BSmc}, Sect.\ 14,   it is shown that there exists
an invariant differential operator
$D$ on $\spX,$ depending on $\tau,$
whose principal symbol is sufficiently generic, such that
\begin{equation}
\naam{e: inversion with D}
D \Wavemc \Foumc = D
\end{equation}
on the Schwartz space $\cC(\spX\col \tau).$
The idea is that $D$ annihilates
the contributions of the discrete and intermediate series to the Plancherel
decomposition of $L^2(\spX\col \tau).$
Accordingly, $\Wavemc\Foumc$ corresponds to the projection onto the most continuous part
of this decomposition at the $K$-type $\tau.$

By an application of {\bf Holmgren's uniqueness theorem} the above mentioned
genericity of the principal symbol of the
differential operator $D$ implies that it is injective on $\Cci(\spX\col \tau),$
see \bib{BSdifop}, Thm.\ 2.
The injectivity of  $\Foumc$ asserted in Theorem \refer{t: injectivity Foumc} follows
from the injectivity of $D$ on $\Cci(\spX\col \tau)$
 combined with (\refer{e: inversion with D}).

In the case of the group Theorem \refer{t: injectivity Foumc} is a straightforward
consequence of the subrepresentation theorem of \bib{CM}.
For indeed, if $f$ belongs to the kernel of $\Foumc,$ then by the subrepresentation theorem,
$f$ is annihilated when integrated against any $K$-finite matrix coefficient.
This is not a valid
argument in the general setting. A priori there might be a $K$-finite
right $H$-fixed generalized matrix coefficient that cannot be produced from the Eisenstein
integrals of the minimal principal series.

We shall now describe the solution to the problem of Fourier inversion mentioned in the introduction.
For this we need the concept of  partial
Eisenstein integral. It follows from the simplifying assumption made
in the beginning of Section \refer{s: Eis}
that
\begin{equation}
\naam{e: defi spXp}
\spXp := K \Aqp H.
\end{equation}
is an open dense subset of $\spX.$
In the situation without the simplifying assumption the definition of $\spXp$
should be adapted by replacing the set on the right-hand side of (\refer{e: defi spXp})
by a finite disjoint union
of open sets of the form $K \Aqp v H,$ with $v$ running through a set $\cW \subset \NKaq$
of representatives for $W/\WKH.$

In the obvious manner we define $\Ci(\spXp\col \tau)$ as the space of
$\tau$-spherical smooth functions
$\spXp \to \Vtau.$ Via restriction, the space (\refer{e: smooth spherical functions}) may naturally be identified
with the subspace of functions in $\Ci(\spXp\col \tau)$ that extend smoothly to the full
space $\spX.$

For $s \in W$ and  $\psi\in \oC$ we define the {\bf partial Eisenstein integral}
$\Eps(\gl\col \dotvar)\psi$ to be the meromorphic function of $\gl \in \faqdc$ with
values in $\Ci(\spXp\col \tau),$ given by
$$
\Eps(\gl\col k a H )\psi  =
 a^{s\gl - \rho}\,\tau(k)\, \Phi_{s \gl} (a)\, [\nC(s \col \gl) \psi](e),\qquad
(a \in \Aqp, k\in K),
$$
for generic $\gl \in \faqdc;$ here $\Phi_\gl$ is as in Proposition \refer{p: asymptotics Eis}.
We agree to write $\Epone = E_{+,1}.$ Then, clearly,
$$
\Eps(\gl \col x)\psi  = \Epone(s\gl\col x) \nC(s \col \gl)\psi, \qquad (\psi \in \oC),
$$
for $s \in W,$ $x\in \spXp$ and generic $\gl \in \faqdc.$
The following result describes the singular set of the functions involved
in the formulation of the inversion theorem. We use the notation (\refer{e: faqd Pmin R}).

\begin{prop}
\naam{p: real singularities}
The functions $\gl \mapsto \dE(\gl\col \dotvar)$ and $\gl \mapsto \Epone(\gl \col\dotvar)$
are meromorphic functions on $\faqdc$
with a singular set consisting of a locally finite union
of hyperplanes of the form $\gl_0 + (\ga^\perp)_\iC,$ with $\gl_0 \in \faqd$ (real)
and with $\ga \in \gS.$ For every $R\in \R$ the set $\faqd(P,R)$
meets only finitely
many of these hyperplanes.
\end{prop}

A proof of this proposition can be found in \bib{BSfi}, Sect.\ 3. Let $\Hyp$
be the collection of singular hyperplanes of $\gl \mapsto \dE(\gl\col \dotvar).$
Then in view of (\refer{e: defi Foumc}) the Fourier transform $\Foumc f$ is meromorphic
on $\faqdc$ with singular locus contained in $\cup \Hyp,$
 for every $f \in \Cci(\spX\col \tau).$

 The solution to the inversion
problem is provided by the following theorem. A sketch of proof 
will be given in Section \refer{s: induction of relations}.

\begin{thm}
\naam{t: Fourier inversion}
{\bf (Fourier inversion theorem){\ }\ } There exists a constant $R < 0$ such that
the functions
$\gl \mapsto \Epone(\gl\col \dotvar)$ and $\gl \mapsto \dE(\gl\col\dotvar)$
are holomorphic in the region $\faqd(\Pmin, R).$
Moreover, let $\eta \in \faqd(\Pmin, R).$ Then, for 
every $f \in \Ci(\spX\col \tau),$
\begin{equation}
\naam{e: Fourier inversion}
f(x) = |W|\; \int_{\eta + i\faqd} \Epone(\gl \col x) \Foumc f(\gl) \; d\gl,
\qquad{\rm for}\quad x \in \spXp.
\end{equation}
\end{thm}

The integral converges absolutely, with local uniformity in $x,$
since the partial Eisenstein integral grows at most of order $(1 + \|\gl\|)^N$ along $\eta + i\faqd,$
for
some $N \in \N,$
whereas the Fourier transform decreases faster than
$C_N (1 + \|\gl\|)^{-N},$
for any $N \in \N.$ Moreover, by Cauchy's integral theorem, the integral
on the right-hand side of (\refer{e: Fourier inversion})
is independent of $\eta$ in the mentioned region. Details can be found in
\bib{BSfi}. The proof of Theorem \refer{t: Fourier inversion}, given in the same paper,
involves  shifting
$\eta$ to $0.$
If no singularities would be encountered during the shift,
then in view of (\refer{e: transformation property Fou})
the integral would become equal to $\Wavemc\Foumc f.$ However,
in general
singularities
are encountered, due to the presence of representations from the discrete and
intermediate series for $\spX.$ This results in residues that can be handled by a calculus that
we developed in \bib{BSres}. These residues can be encoded in terms 
of the concept of
Laurent functional, introduced in the next section. Their contribution to the Fourier inversion
can be analyzed by means of the principle of induction of relations,
also discussed in the next section.

\section{Laurent functionals and induction of relations}
\eqnreset
For the formulation of the principle of induction of relations it is convenient
to introduce the following concept of Laurent functional.

Let $V$ be a finite dimensional real linear space and let $X$ be a finite subset of
$V^*\setminus\{0\}.$ Given a point $a \in V_\iC$ we define the polynomial function
$\pi_a$ on $V_\iC$ by
$$
\pi_a: = \prod_{\xi \in X} (\xi - \xi(a)).
$$
We denote the ring of germs of meromorphic functions 
at $a$ by 
$\cM(V_\iC,a),$ and the subring of germs of holomorphic 
functions by $\cO_a.$ In addition, we define the 
subring
$$ 
\cM(V_\iC, a, X): =  \cup_{N \in \N}\;\;\pi_a^{-N} \cO_a.
$$ 
We now define an $X$-{\bf Laurent functional} at $a \in V_\iC$ 
to be any linear functional
$\Lau \in \Mer(V_\iC, a, X)^*$ such that for every $N \in \N$ there exists a $u_N$ in $S(V),$
the symmetric algebra of $V_\iC,$
such that
\begin{equation}
\naam{e: defi Lau by string}
\Lau = \ev_a \after u_N \after \pi^N_a \qquad {\rm on} \quad \pi_a^{-N} \cO_a.
\end{equation}
Here $S(V)$ is identified with the algebra of translation invariant holomorphic
differential operators on $V_\iC$ and $\ev_a$ denotes evaluation of a function 
at the point $a.$
Finally, the space of all Laurent functionals on $V_\iC,$ relative to $X,$ 
is defined by
\begin{equation}
\naam{e: defi space of laur funct}
\Mer(V_\iC, X)^*_\laur := \bigoplus_{a \in V_\iC} \;\; \Mer(V_\iC, X, a)^*_\laur.
\end{equation}
Given a Laurent functional $\Lau$ from the space on the left-hand side of
(\refer{e: defi space of laur funct}), the finite set of $a \in V_\iC$ for which the
component $\Lau_a$ is non-zero, is called the support of $\Lau$ and denoted by $\supp \Lau.$
Accordingly,
$$
\Lau = \sum_{a \in\, \supp \Lau} \; \Lau_a.
$$
Let now $\Mer(V_\iC,X)$ be the space of meromorphic functions $\gf$
on $V_\iC$ with the property that the germ $\gf_a$ belongs to 
$\Mer(V_\iC,a,X),$ 
for every $a \in V_\iC.$ Then the natural bilinear map 
$(\Lau, \gf) \mapsto \Lau \gf,$ $\Mer(V_\iC,X)^*_\laur \times
\Mer(V_\iC,X) \to \C,$ defined by 
\begin{equation}
\naam{e: Lau on gf} 
\Lau\gf = \sum_{a \in\, \supp \Lau} \; \Lau_a \gf_a,
\end{equation}
induces a linear embedding of $\Mer(V_\iC, X)^*_\laur$ into 
the dual space $\Mer(V_\iC, X)^*.$ 
More details concerning Laurent functionals can be found in \bib{BSanfam}, 
Sect.\ 10.

We end this section with the formulation
of the principle of induction
of relations for the partial Eisenstein integrals $\Eps(\gl\col\dotvar).$
In the proof of the Paley-Wiener theorem, the use of this principle
replaces the use in \bib{Arthur}
of a lifting principle due to W. Casselman, a proof of which
has not appeared in the literature.
Our induction principle does not seem
to imply Casselman's lifting principle for the group. However, it does imply
a version of the lifting principle for normalized Eisenstein integrals,
see \bib{BSanfam}, Thm.\ 16.10.

Let $F\subset \gD$ be a subset of simple roots, let $\gS_F: = \gS \cap \Z F$
be the associated subsystem of $\gS,$ and $W_F$ its Weyl group.
Then $\gS_F$ and $W_F$ are the analogues of $\gS$ and
$W$ for the symmetric space $\spXF = M_F/M_F\cap H.$
Let $W^F \subset W$ be the set of minimal length coset representatives for 
$W/W_F.$
Then the multiplication map of $W$ induces a bijection $W^F\times W_F \to W.$

The group $W_F$ equals the centralizer of $\faFq$ in $W.$
The orthocomplement $\staFq$ of $\faFq$ in $\faq$ is the analogue
of $\faq$ for the space $\spXF.$ Let $K_F = K \cap M_F$ and 
$\tau_F:= \tau|_{K_F}.$
For $t \in W_F$ we denote by
$$
E_{+, t}(\spXF\col \mu \col m) \in \Hom(\oCtau, \Vtau),\qquad 
(\mu \in \staFqdc,\, m \in \spXF),
$$
the analogue for the pair $\spXF , \tau_F$
of the partial Eisenstein integral
$E_{+,t}(\spX\col \gl\col x).$ Here we note that the space $\oCtau$ for $\spX$
coincides with the similar space $\oC(\tau_F)$ for $\spXF.$

\begin{thm}
\naam{t: induction of relations}
{\bf (Induction of relations)\ }
Let, for each $t \in W_F,$ a Laurent functional
$\Lau_t \in \Mer(\staFqdc, \gS_F)^*_\laur \otimes \oC$
be given and assume that
\begin{equation}
\naam{e: relation before induction}
\sum_{t \in W_F} \Lau_t[E_{+, t}(\spXF\col \dotvar\col m)] = 0, \qquad 
(m \in \spXFp).
\end{equation}
Then for each $s \in W^F$ the following meromorphic identity
in the variable $\nu \in \faFqdc$ is valid,
\begin{equation}
\naam{e: relation after induction}
\sum_{t \in W_F}\Lau_t[ E_{+,st}(\spX \col \dotvar + \nu \col x)] = 0, 
\qquad(x \in \spXp).
\end{equation}
Conversely, if (\refer{e: relation after induction}) holds for a fixed $s \in W^F$
and all $\nu$  in a non-empty open subset of $\faFqdc,$
then (\refer{e: relation before induction}) holds.
\end{thm}

This result is proved in our paper \bib{BSanfam}, Thm.\ 16.1.
The proof relies on a more general {\bf vanishing theorem}, see
\bib{BSanfam}, Thm.\ 12.10. This vanishing theorem asserts
that a suitably restricted meromorphic family
$\faFqdc \ni \nu \mapsto f_\nu \in \Ci(\spXp\col \tau)$
of eigenfunctions for $\DX$ is completely determined
by the coefficient of $a^{\nu - \rho_F}$ in its asymptotic expansion
towards infinity along $\AFqp,$ the positive chamber determined by $P_F.$
In particular, if the mentioned
coefficient is zero, then $f_\nu = 0$ for all $\nu;$
whence the name vanishing theorem. Part of the mentioned restriction on
families in the vanishing theorem is a so called {\bf asymptotic globality
condition}. It requires that certain asymptotic coefficients in the 
expansions of $f_\nu$ along
certain codimension one walls should have smooth behavior as functions
in the variables transversal to these walls. The precise condition 
is given in \bib{BSanfam}, Def.\ 9.5.

Let $f^s_\nu,$ for $s \in W^F,$
denote the
expression on the left-hand side of
(\refer{e: relation after induction}).
Then the sum $f_\nu = \sum_{s \in W^F} f^s_\nu$ defines a family
for which the vanishing theorem holds; the summation over $W^F$
is needed for the family to satisfy the asymptotic globality condition.
The
expression on the left-hand side of
(\refer{e: relation before induction}) is the
coefficient of $a^{\nu - \rho_F}$ of the asymptotic 
expansion of $f_\nu$ along $A_{F\iq}^+.$ 
Its vanishing implies that $f = 0.$
From the fact that the sets of the asymptotic exponents of
$f^s_\nu$ along $\AFqp$ are mutually disjoint for distinct $s \in W^F$
and generic $\nu \in \faFqdc,$ it follows that each individual
function $f^s_\nu$ vanishes.
This implies the validity of (\refer{e: relation after induction}).

For the proof of the converse statement it is first
shown that the vanishing of an individual term $f^s_\nu$ implies
that of $f_\nu.$ Here the condition of asymptotic globality once
more plays an essential role. The validity of (\refer{e: relation
before induction}) then follows by taking the coefficient
of $a^{\nu - \rho_F}$ in the asymptotic expansion along $\AFqp.$

\section{Induction of relations and the inversion formula}
\naam{s: induction of relations}
\eqnreset
In this section we shall discuss the role of
induction of relations, as formulated in Theorem
\refer{t: induction of relations}, in the proof of the inversion formula.
Details can be found in \bib{BSfi}.
\medno
{\bf Sketch of proof of Theorem \refer{t: Fourier inversion}{\ }\ }
Let us denote the integral on the right-hand side of (\refer{e: Fourier inversion})
by $\cT_\eta(\Foumc f)(x).$ The main difficulty in the proof is to show
that the function $\cT_\eta \Foumc f \in \Ci(\spXp\col \tau)$ extends
smoothly from $\spXp$ to $\spX.$ By applying a Paley-Wiener shift argument,
with $\eta \to \infty$ in $-\fa_{\iq}^{*+},$
it then follows that $\cT_\eta \Foumc f \in \Cci(\spX\col \tau).$ There exists a differential
operator $D$ as in (\refer{e: inversion with D}),
such that $D \cT_\eta \Foumc f$
is free of singularities during a shift of the integral with $\eta$ moving to $0.$
In view of Cauchy's
theorem this leads to $D\cT_\eta\Foumc f = D \cT_0\Foumc f = D \Wavemc \Foumc f = D f.$ From
the injectivity of $D$ on $\Cci(\spX\col \tau)$ we then obtain (\refer{e: Fourier inversion}).

The most difficult part of the proof concerns the smooth extension of
$\cT_\eta \Foumc f.$ This involves a shift of integration applied to
$\cT_\eta(\Foumc f)(x)$ with  $\eta$ moving to $0.$
According to the residue calculus developed in \bib{BSres},
the process of picking up residues is governed by any choice of a so-called
{\bf residue weight}
on $\gS.$
We fix such a weight, which is by definition a map
$t: \parabs \to [0,1]$ with the property that, for every $Q \in \parabs,$
$$
\sum_{{\buildrel {P \in \parabs} \over { \fa_{P\iq}  = \fa_{Q\iq} } }} t(P) =1.
$$
Moreover, we choose $t$ to be $W$-invariant and even. The latter condition means
that $t(P) = t(\bar P)$ for all $P \in \parabs.$
The encountered residues can be encoded by means of a finite set of Laurent functionals
$\cR_F^t \in \Mer(\faFqdcperp, \Sigma_F)^*_\laur,$ for $F \subset \gD,$
depending only on the root
system $\Sigma,$
the choice of the residue weight $t$ and the locally finite union of hyperplanes which forms
the union of the singular sets of $\gl \mapsto \Epone(\gl\col \dotvar)$
and $\gl \mapsto \dE(\gl\col\dotvar).$

The shift results
in the formula
\begin{eqnarray}\naam{e: cT as residues}
\lefteqn{\cT_\eta (\Foumc f)(x) = }\\\nonumber
&=&
|W|\,\sum_{F \subset \gD} t(P_F)\; \int_{\ge_F + i \faFqd}
\cR^t_F\left(
\sum_{s \in W^F} \Eps(\nu + \dotvar \col x) \Foumc f(\nu + \dotvar) \right) \; d\mu_F(\nu).
\end{eqnarray}
where $\ge_F$ is any choice of elements sufficiently close to zero
in $\faFqdp,$ the positive chamber associated with $P_F.$ Moreover, $d\mu_F$
is the translate by $\ge_F$ of suitably normalized Lebesgue measure on $i \faFqd.$

From the fact that the singular set of the integrand is real in the sense
of Proposition \refer{p: real singularities}, it 
follows that the Laurent functionals
$\cR^t_F$
are real in the following sense. Their support is a set
of real points $a \in \staFqd$ and at each such
point the functional is defined by a string $\{u_N\}\subset S(\staFqd)$
as in (\refer{e: defi Lau by string})
with $u_N$ real for all $N.$

We now define the {\bf kernel functions}
\begin{equation}
\naam{e: formula for the kernel}
K_F(\nu\col x\col y) := \cR^t_F\left( \sum_{s \in W^F}
\Eps(\nu + \dotvar\col x) \dE(\nu + \dotvar \col y) \right).
\end{equation}
Then by using the definition (\refer{e: defi Foumc}) of $\Foumc,$ we may rewrite
the equation (\refer{e: cT as residues}) as
\begin{equation}
\naam{e: deco T eta Foumin f}
\cT_\eta (\Foumc f)(x) = |W|\,\sum_{F \subset \gD}\;t(P_F)\;
 \int_{\ge_F + i \faFqd}
 \left[\int_X K_F^t(\nu \col x\col y) f(y)\; dy \right] \; d\mu_F(\nu).
\end{equation}
For fixed generic $\nu,$
the kernel functions
$K_F^t(\nu \col \dotvar \col \dotvar)
\in \Ci(\spXp \times \spX \col \tau \otimes \tau^*)$  are spherical
and $\DX$-finite in both variables. It follows that they belong
to a tensor product of the form
${\,}^1 E_\nu \otimes {\,}^2 E_\nu,$ with ${}^1 E_\nu$ and ${}^2 E_\nu$
finite dimensional subspaces of $\Ci(\spXp\col \tau)$ and $\Ci(\spXp\col \tau^*),$
respectively. Let ${}^j E_\nu'$ be the subspace of functions
in ${}^j E_\nu$ extending smoothly to $\spX,$ for $j=1,2.$
Then by the symmetry formulated
in Proposition \refer{p: symmetry kernel} below it follows that
the kernel $K_F^t(\nu \col \dotvar\col \dotvar)$
belongs to ${}^1 E_\nu \otimes {\,}^2 E_\nu' \cap {\,}^1 E_\nu' \otimes {\,}^2 E_\nu =
{}^1 E_\nu' \otimes {\,}^2 E_\nu'.$
This shows that the kernel functions extend smoothly to $\spX \times \spX$
and finishes the proof. \qed

\begin{prop}
\naam{p: symmetry kernel}
Let $x,y \in \spXp.$ Then
\begin{equation}
\naam{e: symmetry of the kernel}
K_F^t(\nu \col x\col y) = K_F^t(-\bar \nu \col y\col x)^*
\end{equation}
as a meromorphic identity in the variable $\nu \in \faqdc.$
\end{prop}

Before giving a sketch of the proof we observe that,
due to the fact that $\cR^t_F$ is scalar
and real
in the sense mentioned in the proof of Theorem \refer{t: Fourier inversion} above,
the adjoint of the kernel
is given by
\begin{equation}
\naam{e: adjoint of kernel}
K_F^t(-\bar \nu \col y\col x)^* =
\cR^t_F\left( \sum_{s \in W^F}
\nE(\nu - \dotvar\col x) \dEps(\nu - \dotvar \col y) \right),
\end{equation}
where the dual partial Eisenstein integrals are defined by
$$
\dEps(\gl \col x) := \Eps(- \bar \gl \col x)^*.
$$

{\bf Sketch of proof of Proposition \refer{p: symmetry kernel}{\ }\ }
The final part of the proof of Theorem \refer{t: Fourier inversion}
can be modified in such a way that (\refer{e: symmetry of the kernel})
is only needed for $F \subset \gD$ with $F \neq \gD.$ The validity
of (\refer{e: symmetry of the kernel}) for $F = \gD$ is derived in the course
of the modified argument. For details, we refer the reader to
\bib{BSfi}, Sect.\ 9.

Thus,
we may restrict ourselves to proving (\refer{e: symmetry of the kernel})
for $F \subsetneq \gD.$ This in turn is achieved by using induction
of relations in order to reduce to the lower dimensional
space $\spX_F.$ More precisely, the residue weight $t$ naturally induces
a residue weight ${}^*t$ on $\gS_F,$ the analogue of $\gS$ for $\spX_F.$
The set $F$ is a simple system for $\gS_F.$ Let
$K_F^{\start}(\spXF\col \dotvar \col \dotvar)$
be the analogue of $K^t_\gD$ for the space $\spX_F.$
Then by induction, $K_F^{\start}(\spXF\col \dotvar \col \dotvar)$
is a smooth function on
$\spX_F \times \spX_F$ and satisfies the symmetry condition
\begin{equation}
\naam{e: symmetry of the lower rank kernel}
K^{\start}_F(\spXF\col m \col m') = K^{\start}_F(\spXF \col m' \col m)^*,
\end{equation}
for $m,m' \in \spXF.$
Here we have suppressed the analogue of the parameter $\nu,$ which
is zero dimensional in the present setting.

The residue calculus behaves well with respect to induction.
In particular, let $\cR_F^{\start} \in \Mer(\staFqd, \gS_F)^*_\laur$
be the analogue of $\cR_\gD^t$ for the data $\spXF, \gS_F, F, \start.$
Then $\cR_F^{\start} = \cR_F^t;$ for obvious reasons,
we have called this
result {\bf transivity of residues}, see \bib{BSres}, Sect.\ 3.6.
Using (\refer{e: formula for the kernel}) and (\refer{e: adjoint of kernel})
for $K^{\start}_F(\spXF),$ taking into account that 
$(W_F)^F = \{1\},$ we thus see that 
(\refer{e: symmetry of the lower rank kernel})
is equivalent to
\begin{eqnarray}\nonumber\lefteqn{
\cR^t_F\left(
E_{+}(\spXF\col \dotvar\col m) \dE(\spXF \col \dotvar \col m') \right)
}\qquad\qquad\\  &=&
\cR^t_F\left(
\nE(\spXF \col - \dotvar \col m) \dE_+(\spXF\col - \dotvar\col m')  \right),
\naam{e: symmetry kernel written out}
\end{eqnarray}
where $\dE_+ := \dE_{+,1}.$
In view of (\refer{e: formula for the kernel}) and (\refer{e: adjoint of kernel}), the relation
(\refer{e: symmetry of the kernel}) can now be derived from
(\refer{e: symmetry kernel written out}), by applying induction
of relations, first with respect to the variable $x$ and then
a second time
with respect to the variable $y.$ For details we
refer the reader to \bib{BSfi}, Sect.\ 8.
\qed

\section{Arthur-Campoli relations}
\eqnreset
In this section we describe the so called Arthur-Campoli relations,
needed for the formulation of the Paley-Wiener theorem in the
next section. We start with the definition of an {Arthur-Campoli functional}.

\begin{defi}
An {\bf Arthur-Campoli functional} for $\spX, \tau$
is a Laurent functional $\Lau \in \Mer(\faqdc, \Sigma)^*_\laur \otimes \oCtau$ with the property that
$$
\Lau \dE( \dotvar \col x) = 0 \text{for all} x \in \spX.
$$
The linear space of such functionals is denoted by $\AC(\spX\col \tau).$
\end{defi}

From the principle of induction of relations
as formulated in Theorem \refer{t: induction of relations},
the following result follows in a straightforward manner.
See \bib{BSpw} for details.

\begin{lemma}
\naam{l: induction of AC rels}
{\bf (Induction of AC relations)\ }{\ }
Let $F \subset \gD$ and $\Lau \in \AC(\spXF\col \tau_F).$
Then for generic $\nu \in \faFqdc,$ the Laurent functional
$$
\Lau_\nu: \gf \mapsto \Lau[\gf(\nu + \dotvar)]
$$
belongs to $\AC(\spX\col \tau).$
\end{lemma}

In this result, `generic' can  be made more precise as follows. There exists
a locally finite collection $\Hyp_S$ of hyperplanes in $\faFqdc,$ specified explicitly in terms of
the support $S$ of $\Lau,$ such that the statement is valid for
$\nu \in \faFqdc \setminus \cup \Hyp_S.$

\section{The Paley-Wiener theorem}
\eqnreset
In this section we shall formulate the Paley-Wiener
theorem, and indicate how induction of relations enters
its proof.
Our first objective is to define a space of Paley-Wiener functions.
The first step is to define
a suitable space of meromorphic functions
that takes the singularities of the Fourier transform into account.

Let $\cH = \cH(X,\tau)$ be the smallest collection of hyperplanes of
the form $\gl_0 + (\ga^\perp)_\C,$ with $\gl_0 \in \faqd$ and $\ga \in \gS,$
such that the $\Ci(\spX) \otimes \Hom(\Vtau, \oCtau)$-valued  meromorphic function
$\gl \mapsto \dE(\gl\col \dotvar)$ is regular on $\faqdc \setminus \cup \cH.$
By the requirement of minimality,
the collection $\cH$ has the properties of Proposition \refer{p: real singularities}.

If $H \in \cH$ we select $\ga_H \in \gS$ and $s_H \in \R$ such that
$H$ is given by the equation $\inp{\gl}{\ga_H} = s_H.$
Let $d(H)$  denote the order of the singularity of
$\gl \mapsto \dE(\gl)$ along $H.$ Thus, $d(H)$ is the smallest natural
number for which $\gl \mapsto (\inp{\gl}{\ga_H} - s_H)^{d(H)} \dE(\gl)$
is regular at the points of $H$ that are not contained in any hyperplane
from $\Hyp \setminus \{H\}.$

If $\omega \subset \faqdc$ is a bounded subset, 
then in view of the mentioned properties of $\Hyp$ 
we may define a polynomial function $\pi_\omega: \faqdc \to \C$ by
$$
\pi_\omega(\gl) = \prod_{H \in \Hyp, H \cap \omega \neq \emptyset}\;\;
(\inp{\gl}{\ga_H} - s_H)^{d(H)}.
$$
We define $\Mer(\faqdc, \Hyp, d)$ to be the space of meromorphic functions
$\gf: \faqdc \to \C$ such that, for every bounded open set $\omega \subset \faqdc,$
the function $\pi_\omega \gf$ is regular on $\omega.$ Taking into account
that the $\ga_H$ and $s_H$ are real for $H \in \Hyp,$ we readily see
that for each function $\gf \in \Mer(\faqdc, \Hyp, d)$ and every bounded
open subset $\omega \subset \faqdc,$ the function $\pi_\omega \gf$
is in fact regular on $\omega + i \faqd.$

In view of the definitions just given, the function $\gl \mapsto \dE(\gl\col x)$
belongs to the space $\Mer(\faqdc, \Hyp, d) \otimes \Hom(\Vtau, \oCtau),$
for every $x \in \spX.$ Moreover, $\Foumc$ maps $\Cci(\spX\col \tau)$ into
$\Mer(\faqdc, \Hyp, d) \otimes \oCtau.$

It follows from Proposition \refer{p: real singularities}
that the set $\Hyp_0$ of $H \in \Hyp$ having empty intersection with
${\rm cl}\,\faqd(\Pmin,0)$ is finite. We define
the polynomial function $\pi: \faqdc \to \C$ by
$$
\pi(\gl) = \prod_{H \in \Hyp_0} ( \inp{\gl}{\ga_H} - s_H)^{d(H)}.
$$
Then there exists a constant $\ge >0$ such that $\gl \mapsto \pi(\gl) \dE(\gl)$
is regular on $\faqd(\Pmin, \ge).$ It follows that for every $f \in \Cci(\spX\col \tau)$
the $\oCtau$-valued meromorphic function $\gl \mapsto \pi(\gl) \Foumc f(\gl)$
is regular on $\faqd(\Pmin, \ge).$

We now define $\cP(\faqdc, \Hyp, d)$ as the subspace of $\Mer(\faqdc, \Hyp, d)$ consisting
of functions $\gf$ which satisfy the following condition of decay in the imaginary directions
$$
\sup_{\gl \in \omega + i \faqd}\; (1+ |\gl|)^n \,| \pi_\omega(\gl) \gf(\gl)| < \infty,
$$
for every compact set $\omega \subset\faqd$ and all $n \in \N.$
Equipped with the suggested seminorms, the space $\cP(\faqdc, \Hyp, d)$ is a Fr\'echet
space. Moreover, via (\refer{e: Lau on gf})
the space  of Laurent functionals $\Mer(\faqdc, \gS)^*_\laur$
naturally embeds into the continuous linear dual of $\cP(\faqdc, \Hyp, d).$
It follows that the following subspace of $\cP(\faqdc, \Hyp, d) \otimes \oCtau$
is closed, hence Fr\'echet,
$$
\cP_\AC(\spX\col \tau): = \{\gf \in \cP(\faqdc, \Hyp, d) \otimes \oCtau
\mid
\Lau \gf = 0,\;\; \forall \; \Lau \in \AC(\spX\col \tau)\}.
$$
Finally, we define the Paley-Wiener space by incorporating a condition of exponential growth
along a closed cone.

\begin{defi}
The {\bf Paley-Wiener space} $\PW(\spX\col \tau)$
is defined to be the space of functions
$\gf \in \cP_\AC(\spX \col \tau)$ for which there exists a constant $M > 0$
such that, for all $n \in \N,$
$$
\sup_{\gl \in \,{\rm cl}\,\faqd(\Pmin, 0)}\; ( 1 + |\gl|)^n \,e^{-M |\Re \gl|}\, \|\pi(\gl) \gf(\gl)\| < \infty.
$$
The subspace of functions satisfying this estimate with a fixed $M > 0$ and all $n \in \N$
is denoted by $\PWM(\spX\col \tau).$
\end{defi}

By using Euclidean Fourier analysis, it can be shown that
$\PWM(\spX\col \tau)$ is a closed subspace of $\cP_\AC(\spX\col \tau),$
for each $M >0,$ hence a Fr\'echet space for the restriction topology.
For details we refer the reader to \bib{BSpw}.
Accordingly, for $M < M'$ we have a continuous linear embedding
of $\PWM(\spX\col \tau)$ onto a closed subspace of ${\rm PW}_{M'}(\spX\col \tau).$
The space $\PW(\spX\col \tau),$ being the union of the spaces $\PWM(\spX\col \tau),$
is equipped with the associated direct limit topology. Thus, it becomes a
strict {\rm LF}-space.

For $M >0$ we denote by $B_M$ the closed ball in $\faq$ of center $0$ and radius $M.$
Moreover, we denote by $\CiM(\spX\col \tau)$ the space of functions
in $\Ci(\spX\col \tau)$ with compact support contained in $K \exp B_M H.$

\begin{thm}{\bf (Paley-Wiener theorem)\ }{\ }
\naam{t: PW}
The Fourier transform $\Foumc$ is a topological linear
isomorphism from $\Cci(\spX\col \tau)$ onto $\PW(\spX\col \tau).$ More precisely,
for each $M >0$ it maps $\CiM(\spX\col \tau)$ homeomorphically onto
$\PWM(\spX\col \tau).$
\end{thm}

In the Riemannian case $H = K$ and $\tau = 1,$ this result is equivalent to the Paley-Wiener
theorem of S.\ Helgason and R.\ Gangolli, see \bib{Hel}, Thm.\ IV, 7.1.
In the case of the group our Paley-Wiener theorem can be shown to be equivalent
to the one of J.\ Arthur, \bib{Arthur}, which in turn generalizes the result of O.A.\ Campoli,
\bib{Campoli},
for groups of split rank one. Arthur's proof relies on Harish-Chandra's
Plancherel theorem and the lifting principle mentioned in
Section \refer{s: induction of relations}, due to W.\ Casselman. It also makes
use of ideas from the residue calculus appearing 
in the work of R.P.\ Langlands, \bib{LEis}.
In \bib{DPW}, P.\ Delorme used a different method to obtain a Paley-Wiener 
theorem for semisimple groups with one conjugacy class of Cartan subgroups, 
with explicit symmetry conditions instead of the Arthur-Campoli relations. 
This work
in turn generalized work of Zhelobenko, \bib{Zhe}, for the complex groups.

We conjectured the present Paley-Wiener
theorem in slightly different but equivalent form in \bib{BSmc}, where we proved
it under the assumption that $\dim \faq =1.$ 
The proof of Theorem \refer{t: PW}
is given in the paper \bib{BSpw}. It relies on the inversion
theorem, Theorem \refer{t: Fourier inversion}, and on the principle of induction
of relations, see Theorem \refer{t: induction of relations}. In particular, our
proof is independent of the theory of the discrete series
and the existing proofs of the Plancherel theorem (in \bib{Dpl}, \bib{BSpl1} and \bib{BSpl2}).
The precise
relation with the Plancherel decomposition will be described in Section \refer{s: Plancherel}.

In the following
sketch we will indicate the main ideas of our proof of the Paley-Wiener theorem.
\medbreak
{\bf Sketch of proof of Theorem \refer{t: PW}\ }{\ }
As usual, the proof that $\Foumc$ maps $\CiM(\spX\col \tau)$ continuously into
$\PWM(\spX\col \tau)$ is rather straightforward. For details, see
\bib{BSmc}. The injectivity of $\Foumc$ was already
asserted in Theorem \refer{t: injectivity Foumc}.
By the open mapping theorem for Fr\'echet spaces, it remains to establish
the surjectivity of $\Foumc.$ Let $\gf \in \PWM(\spX\col\tau).$
In view of the inversion theorem the only possible candidate
for a function $f \in \CiM(\spX\col \tau)$
with Fourier transform equal to $\gf$ is given by the formula
$$
f(x) =
|W|\,\int_{\eta + i\faqd} \Epone(\gl \col x) \gf(\gl) \; d\gl,
$$
for $x \in \spXp$ and for $\eta \in \faqd$ sufficiently $\bar\Pmin$-dominant.
The problem with this formula is that it only defines a smooth function $f$
on the open dense subset $\spXp$ of $\spX.$ By a standard shift argument
of Paley-Wiener type,
with $\eta$ moving to infinity in $- \faqdp,$ it follows that the support
of $f$ is contained in $K \exp B_M H.$ Therefore, it suffices to show that
the function $f$ has a smooth extension to all of $\spXp.$ This is the central
theme of the proof.

We will actually show that $f$ has a smooth extension under the weaker assumption
that $\gf \in \cP_\AC(\spX\col \tau).$
As in the proof of Theorem \refer{t: Fourier inversion} the idea is to write
the integral differently by application of a contour shift, with $\eta$ moving to $0,$
and by organizing the residual integrals according to the calculus described in the mentioned
proof. This leads to the formula
\begin{equation}
\naam{e: inversion formula on PW space}
f(x) = \sum_{F\subset \gD}  \cT_F^t \gf (x), \qquad (x \in \spXp),
\end{equation}
with
\begin{equation}
\naam{e: definition cT F}
\cT_F^t \gf(x): = |W|\, t(P_F) \int_{\ge_F + i \faFqd}
\cR^t_F\left(
\sum_{s \in W^F} \Eps(\nu + \dotvar \col x) \gf(\nu + \dotvar) \right) \; d\mu_F(\nu).
\end{equation}
The problem now is to show that each of the individual terms $\cT_F \gf$
extends smoothly to all of $\spX.$ This is done by writing $\cT_F \gf$
as a superposition of certain generalized Eisenstein integrals.

These were defined in \bib{BSfi} by using the symmetry property
of the kernels $K_F^t,$ as formulated in Proposition \refer{p: symmetry kernel}.
As in the proof of Theorem \refer{t: Fourier inversion} let $K_F^{\start}(\spXF)
 \in \Ci(\spXF \times \spXF) \otimes
\End(\Vtau)$ be the analogue for $\spXF$ and  $\tau_F$
 of the kernel $K_\gD^t$ for $\spX$
and $\tau.$ We recall that $K_F^{\start}(\spXF)$ does not depend
on a spectral parameter, since the analogue of $\fa_{\gD\iq}$ for $\spXF$ is the zero space.
We define the following subspace of $\Ci(\spXF\col \tau_F),$
$$
\cA_F = \cA^{\start}(\spXF\col \tau_F): =
{\rm span}\, \{K_F^{\start}(\spXF \col \dotvar \col m')u \mid m' \in \spXFp,
u \in \Vtau\}.
$$
Being annihilated by a cofinite ideal of $\D(\spXF),$ this space is finite dimensional.
It can be shown that $\cA_F$ is the discrete
series subspace $L^2_d(\spXF\col \tau_F)$ of $L^2(\spXF\col \tau_F),$  see \bib{BSpl1},
Lemma 12.6 and Thm.\ 21.2,
but this fact is not needed
for the proof of the Paley-Wiener theorem.

For $\psi \in \cA_F$ we define the {\bf generalized Eisenstein integral}
$\nE_F(\psi \col \nu)$ as a meromorphic $\Ci(\spX\col \tau)$-valued function
of $\nu \in \faFqdc,$ as follows.
If
\begin{equation}
\naam{e: representation of psi in cA F}
\psi = \sum_i K^\start_F(\spXF\col \dotvar \col m_i') u_i,
\end{equation}
with $m_i' \in \spXFp$ and $u_i \in \Vtau,$ then
\begin{equation}
\naam{e: defi nEF}
\nE_F(\psi\col \nu \col x): =
\sum_i \cR^t_F [ \nE(\nu - \dotvar \col x) E_+^*(\spXF \col - \dotvar \col m_i')u_i].
\end{equation}
It follows by induction of relations, Theorem \refer{t: induction of relations},
 that the expression
(\refer{e: defi nEF}) is independent of the particular representation
of $\psi \in \cA_F$ given in (\refer{e: representation of psi in cA F}).
It also follows by induction of relations, combined with the symmetry of the
kernel $K^\start_F,$ that for $\psi \in \cA_F$ given by (\refer{e: representation of psi in cA F}),
\begin{equation}
\naam{e: alternative representation nE F}
\nE_F(\psi \col \nu \col x) = \sum_i \cR_F^t [
\sum_{s \in W^F} E_{+,s}(\nu + \dotvar \col x)E^*(\spXF \col\dotvar \col m_i')u_i],
\end{equation}
for generic $\nu \in \faFqdc$ and all $x \in \spXp.$
Let 
$$
T_F(\spXF\col \dotvar) = T_F^\start(\spXF\col \dotvar)\,:\;
 \Cci(\spXF\col \tau_F) \to \cA_F
$$ 
be the analogue for $\spXF$
of
the operator $T_\gD^t$ occurring in (\refer{e: inversion formula on PW space}).
Then it follows from (\refer{e: defi nEF}) and (\refer{e: alternative representation nE F}),
essentially by integration with respect to the variable $m'$ that,
for all $f \in \Cci(\spXF\col \tau_F),$
\begin{equation}
\naam{e: nEF of TF}
|W_F|^{-1}\; \nE_F( T_F(\spXF\col f) \col \nu \col x)
= \cR_F^t[ \sum_{s \in W^F} E_{+, s}(\nu + \dotvar \col x) \Foumc(\spXF \col f)(\dotvar)].
\end{equation}
Here $\Foumc(\spXF\col \dotvar)$ denotes the analogue of $\Foumc$ 
for $\spXF.$ 

The next step in the proof of the Paley-Wiener theorem consists of the following
result, which follows from the Arthur-Campoli relations and their inductive
property described in Lemma \refer{l: induction of AC rels},
essentially by application of linear algebra.

\begin{prop}
\naam{p: application of AC to res int}
Let $F \subset \gD.$ There exists a finite dimensional complex linear subspace
$V \subset \Cci(\spXF\col \tau_F)$ and a Laurent functional
$\Lau' \in \Mer(\staFqdc, \gS_F)^*_\laur \otimes \Hom(\oCtau, V)$
such that, for generic $\nu \in \faFqdc,$ the map
$
\gf \mapsto f_{\nu, \gf},$
$\cP_\AC(\spX\col \tau) \to V,$  defined by
$$
f_{\nu, \gf} = \Lau'[\gf(\nu + \dotvar)],
$$
has the following property, for all $x \in \spXp,$
\begin{eqnarray*}
\lefteqn{
\cR^t_F \left[ \sum_{s \in W^F} E_{+,s}(\nu + \dotvar\col x) 
\gf(\nu + \dotvar)\right]}\qquad\qquad\qquad
\\
\qquad\qquad\qquad&=&
\cR^t_F \left[ \sum_{s \in W^F} E_{+,s}(\nu + \dotvar\col x)
\Foumc(\spXF\col f_{\nu, \gf})(\dotvar)\right].
\end{eqnarray*}

\end{prop}

The final step in the proof is the following result, which follows
by combining Proposition \refer{p: application of AC to res int}
with (\refer{e: nEF of TF}).

\begin{prop}
There exists a  $\Lau_F \in \Mer(\staFqdc, \gS_F)^*_\laur\otimes \Hom(\oCtau, \cA_F)$
such that
$$
\cR^t_F\left[ \sum_{s \in W^F} E_{+,s}(\nu + \dotvar\col x)\gf(\nu + \dotvar)\right]
=
\nE_F(\Lau_F[\gf(\nu + \dotvar)]\col \nu \col x),
$$
for all $\gf \in \cP_\AC(\spX\col \tau),$ $x \in \spXp$ and generic $\nu \in \faFqdc.$
\end{prop}

It follows from combining this proposition with (\refer{e: definition cT F})
that, for $\gf \in \cP_\AC(\spX\col \tau),$
\begin{equation}
\naam{e: inversion formula per F on P AC space}
\cT_F \gf (x) = |W|\, t(P_F)\,
\int_{\ge_F + i \faFqd} \nE_F(\Lau_F[\gf(\nu + \dotvar)]\col \nu \col x) \; d\mu_F(\nu),
\end{equation}
for all $x \in \spXp.$ From this expression it is readily seen that
$\cT_F$ extends to a continuous linear map $\cP_\AC(\spX\col \tau) \to \Ci(\spX\col \tau).$
\qed

\section{Relation with the Plancherel decomposition}
\naam{s: Plancherel}
\eqnreset
In this section we briefly discuss the relation between the Paley-Wiener theorem
and the Plancherel theorem, obtained by P.\ Delorme \bib{Dpl} and, independently,
by H.\ Schlichtkrull and myself in \bib{BSpl1} and \bib{BSpl2}. Earlier, a
Plancherel theorem had been announced by T. Oshima, \bib{Opl}, p.\ 32,
but the details have not appeared.
For the case of the group, the Plancherel theorem is due to Harish-Chandra,
\bib{HC1}, \bib{HC2}, \bib{HC3}.
For the case of a complex reductive group modulo
a real form, the Plancherel theorem has been obtained by P.\ Harinck, \bib{Hpl}.

The starting point of our proof of the Plancherel theorem is the Fourier inversion
formula
\begin{equation}
\naam{e: inversion formula with kernels}
f(x) = |W|\, \sum_{F \subset \gD} t(P_F) \int_{\ge_F + i \faFqd} \,\int_\spX
K^t_F(\nu \col x \col y) f(y) \; dy \,d\mu_F(\nu), \qquad (x \in \spX),
\end{equation}
which follows from Theorem \refer{t: Fourier inversion} and (\refer{e: deco T eta Foumin f}).
The crucial part of the proof of the Plancherel theorem consists of showing that this formula,
which is valid for $\ge_F$ sufficiently close to zero in $\faFqdp,$ remains valid
with $\ge_F = 0$ for all $F \subset \gD.$ This in turn is achieved by
showing that the kernel functions
$K^t_F$ are regular for $\nu \in i\faFqd.$

The regularity is achieved in a long inductive
argument in \bib{BSpl1}. It is in this argument that we need the theory
of the {\bf discrete series} for $\spX$
initiated by M.\ Flensted-Jensen \bib{FJds} and further developed
in the fundamental paper \bib{OMds} by T.\ Oshima and T.\ Matsuki. Of the latter paper
two results on the discrete series are indispensable.
The crucial results needed are the necessity and sufficiency of the rank condition for the discrete
series to be non-empty as well as the fact that  representations from the discrete series
have real and regular $\DX$-characters; see \bib{BSpl1} for details.

In the course of the inductive argument, it is is shown that $K_F^t$ is independent of
the choice of the residue weight $t;$ moreover, $\cA_F = L^2_d(\spXF\col \tau_F)$
and the generalized Eisenstein integral $\nE_F$ is independent of $t$ as well.
It is then shown that
\begin{equation}
\naam{e: kernel as prod Eis}
K_F(\nu \col x \col y) = |W_F|^{-1} \nE_F(\nu \col x)E_F^*( \nu \col y),
\end{equation}
with $E_F^*( \nu \col y):= \nE_F(-\bar \nu \col y)^*.$
At this point we note that if we define the Fourier transform
$\cF_F: \Cci(\spX\col \tau_F) \to \Mer(\faFqdc) \otimes \cA_F$
as $\Foumc$ in (\refer{e: defi Foumc}) with $\dE_F$ in place of $\dE,$ then
(\refer{e: inversion formula with kernels}) becomes
\begin{equation}
\naam{e: Fourier inversion with Eis}
f(x) = \sum_{F \subset \gD} [W\col W_F] \,t(P_F) \int_{\ge_F + i\faFqd} 
\nE_F(\Fou_F f(\nu) \col \nu \col x)\; d\mu_F(\nu).
\end{equation}
The relation of this formula with (\refer{e: inversion formula on PW space}) and
(\refer{e: inversion formula per F on P AC space})
for $\gf = \Foumc f$ is given by
$$
\nE_F(\Fou_F f(\nu) \col \nu \col x) = |W_F|^{-1}\,\nE_F(\Lau_F[\Foumc f] (\nu + \dotvar)\col \nu \col x),
$$
for every $x \in \spX,$ as an identity of meromorphic functions in the variable
 $\nu \in \faFqdc.$
From taking coefficients of $a^{\nu - \rho_F}$ in the asymptotic expansions of both members
along $M_F \AFqp$
it follows that
$$
\Fou_F f(\nu) = |W_F|^{-1}\,  \Lau_F[\Foumc f](\nu + \dotvar),
$$ for every $f \in \Cci(\spX\col \tau),$
as an identity of meromorphic functions in the variable $\nu \in \faFqdc.$ This in turn leads to
the meromorphic identity
$E^*_F(\nu \col x) = |W_F|^{-1}\, \Lau_F[\dE(\nu + \dotvar\col x)],$ for all $x \in \spX.$

In view of (\refer{e: kernel as prod Eis}),
the regularity result for the kernel is reduced to the similar result for
the generalized Eisenstein integral $\nE_F(\nu \col \dotvar),$ namely its regularity
for $\nu \in i\faFqd.$ This is the analogue of Theorem \refer{t: regularity}.
By the work of J.\ Carmona on the theory of the {\bf constant term} for $\spX,$ which
in turn generalizes Harish-Chandra's work \bib{HC1} for the case of the group,
we can define generalized $c$-functions, which are the analogues of the $c$-functions
in Proposition \refer{p: asymptotics Eis}. A key step in the proof of the regularity theorem
 is then to prove
the {\bf Maass-Selberg relations} for these generalized $c$-functions, see Theorem
\refer{t: Maass-Selberg}. It should be said that at the time of the announcement
of our proof of the Plancherel theorem we had to rely on the  Maass-Selberg
relations proved by Delorme in \bib{Dtr}. Since then we have found a way to derive
the generalized Maass-Selberg relations from those associated with a minimal
$\gs\Cartan$-stable parabolic subgroup, as formulated in Theorem \refer{t: Maass-Selberg};
see \bib{BSpl1}, Thm. 18.3.

From the regularity theorem it follows that (\refer{e: Fourier inversion with Eis})
holds with $\ge_F = 0.$
Defining the wave packet transform $\Wave_F$ as $\Wavemc$
in (\refer{e: defi wave packet}) with $\nE_F$ instead of $\nE$ we now obtain that
\begin{equation}
\naam{e: inversion with ge F zero}
f = \sum_{F \subset \gD} [W\col W_F]\, t(P_F) \Wave_F \Fou_Ff.
\end{equation}
In \bib{BSpl1} we establish uniform tempered estimates for the generalized Eisenstein integrals.
These allow to show that the formula (\refer{e: inversion with ge F zero})
extends continuously to the Schwartz space
$\cC(\spX\col \tau).$  It can be shown that $\Wave_F\after \Fou_F$
depends on $F$ through its class for the
 equivalence relation $\sim$ on the powerset $2^\gD$ defined by
$F \sim F'\iff \exists w \in W: \, w(\faFq) = \fa_{F'\iq}.$
By a simple counting argument it then follows that
\begin{equation}
\naam{e: Plancherel for spherical functions}
I = \sum_{[F] \in 2^\gD/\sim} [W\col W^*_F]\; \Wave_F\Fou_F \qquad \text{on} \cC(\spX\col \tau);
\end{equation}
here $W^*_F$ denotes the normalizer of $\faFq$ in $W.$ In particular,
in this {\bf Plancherel formula} for $\tau$-spherical functions the residue
weight $t$ has disappeared.

In \bib{BSpl2} it is shown that the Eisenstein integrals $\nE_F(\nu),$ 
for $\nu \in \faFqdc,$
are essentially
sums of generalized matrix coefficients of parabolically induced representations of the form
$\Ind_{P_F}^G(\gs\otimes \nu \otimes 1)$
with $\gs$ a discrete series representation of $\spXF = M_F/M_F \cap H.$
Here a key role is played by the {\bf automatic continuity
theorem} due to W.\ Casselman and N.\ Wallach, \bib{Cas} and \bib{Wal2}.
This allows to conclude that (\refer{e: Plancherel for spherical functions})
is the $\tau$-spherical part of the Plancherel formula in the sense of representation
theory. Moreover, the Eisenstein integrals $\nE_F(\nu)$ and the associated
Fourier and wave packet transforms can be identified with those
introduced in \bib{CDn} by Carmona and Delorme.
%
\def\adritem#1{\hbox{\small #1}}
\def\distance{\hbox{\hspace{8.5cm}}}
\def\apetail{@}

\def\banaddress{\vbox{
\adritem{E.\ P.\ van den Ban}
\adritem{Mathematisch Instituut}
\adritem{Universiteit Utrecht}
\adritem{PO Box 80 010}
\adritem{3508 TA Utrecht}
\adritem{Netherlands}
\adritem{E-mail: ban{\apetail}math.uu.nl}
}
}

\vbox{\vspace{1mm}}
\vfill
\hbox{\vbox{\distance}\hfill\vbox{\banaddress}}
\end{document}